\documentclass[article]{ij4uq_arxiv}%
\usepackage{amsmath}
\usepackage{amsfonts}
\usepackage{amssymb}
\usepackage{graphicx}
\usepackage{moreverb}
\usepackage{multirow}
\usepackage{subfigure}
\usepackage{url}
\usepackage[pdftex]{thumbpdf}
\usepackage[pdftex,
pdfstartview=FitH,bookmarks=true,bookmarksnumbered=true,bookmarksopen=true,hypertexnames=false,breaklinks=true,
colorlinks=true,  linkcolor=blue,anchorcolor=blue,
citecolor=blue,filecolor=blue,  menucolor=blue]{hyperref}%
\providecommand{\U}[1]{\protect\rule{.1in}{.1in}}
\providecommand{\U}[1]{\protect\rule{.1in}{.1in}}
\newtheorem{theorem}{Theorem}
\newtheorem{algorithm}[theorem]{Algorithm}

\fancypagestyle{plain}{
\fancyhf{}
\fancyfoot[R]{\vspace*{10pt}\small\bf\thepage}
\fancyfoot[L]{\fottitle}
}
\begin{document}
\volume{Volume 1, Number 1, \myyear\today}
\title{Truncated hierarchical preconditioning for the stochastic {G}%
alerkin {FEM}}
\titlehead{Truncated Hierarchical Preconditioning}
\authorhead{B. Soused\'{\i}k \& R. G. Ghanem}
\corrauthor[1,2]{Bed\v{r}ich Soused\'{i}k}
\corremail{\href{mailto:sousedik@usc.edu}{sousedik@usc.edu}}
\author[1]{Roger G. Ghanem}
\email{ghanem@usc.edu}
\corrurl{\url{http://venus.usc.edu/}}
\address[1]{Department of Aerospace and Mechanical Engineering \\
University of Southern California  \\
Los~Angeles, CA 90089-2531, USA}
\address[2]{Institute of Thermomechanics \\
Academy of Sciences of the Czech Republic \\
Dolej\v{s}kova 1402/5, 182~00 Prague~8, Czech Republic}
\dataO{}
\dataF{}
\abstract{
Stochastic Galerkin finite element discretizations of partial differential equations
with coefficients characterized by arbitrary distributions lead, in general, to fully block dense linear systems.
We propose two novel strategies for constructing preconditioners for these systems
to be used with Krylov subspace iterative solvers.
In particular, we present a variation on of the hierarchical Schur complement preconditioner, developed recently by the authors,
and an adaptation of the symmetric block Gauss-Seidel method.
Both preconditioners take advantage of the hierarchical structure of global stochastic Galerkin matrices,
and also, when applicable, of the decay of the norms of the stiffness matrices obtained from the polynomial chaos expansion of the coefficients.
This decay allows to truncate the matrix-vector multiplications in the action of the preconditioners.
Also, throughout the global matrix hierarchy, we approximate solves with certain submatrices by the associated diagonal block solves.
The preconditioners thus require only a limited number of stiffness matrices
obtained from the polynomial chaos expansion of the coefficients,
and a preconditioner for the diagonal blocks of the global matrix.
The performance is illustrated by numerical experiments.
}
\keywords{
stochastic Galerkin finite element methods,
iterative methods,
Schur complement method,
Gauss-Seidel method,
hierarchical and multilevel preconditioning
}
\maketitle

\section{Introduction}

\label{sec:introduction}Precise values of coefficients in the setup of
physical models using partial differential equations (PDEs) are often not
known. In such situations, the coefficients are typically treated as random
variables or processes in attempt to quantify uncertainty in the underlying
problem. Probably the most popular and widely used class of methods for a
solution of these problems are the Monte-Carlo techniques, because they
require only solutions of the PDE for a given set of realizations of the input
random coefficients. These methods are well-known for their robustness,
versatility, and quite slow convergence. Therefore, in the last two decades, a
significant effort has been devoted to the development of methods that
leverage regularity of the solution,\ and outperform the Monte-Carlo methods
at least for problems with stochastic dimensions that are not too large. The
most promising developments include the stochastic finite element methods.
There are two main variants of the stochastic finite elements: collocation
methods~\cite{Babuska-2010-SCM,Xiu-2005-HCM}, and stochastic Galerkin
methods~\cite{Babuska-2005-GFE,Babuska-2005-SEB,Ghanem-1991-SFE}. The first
approach samples the stochastic PDE at a predetermined set of collocation
points, which yields a set of uncoupled deterministic systems. The solution at
these collocation points is then used to interpolate the solution in the
entire random input domain. Because extending legacy software to support the
collocation points is relatively simple, these methods are often regarded as
non-intrusive. On the other hand, the second approach, using the stochastic
Galerkin methods, translates the stochastic PDE into one large coupled
deterministic system. Consequently, the Galerkin methods require a development
of new solvers and therefore they are commonly regarded as intrusive. Since
use of direct solvers for the stochastic Galerkin systems might be prohibitive
due to their large size, iterative solvers are often preferred. We then seek
preconditioners to speed up convergence.

There are two main types of expansion of the random field representing the
uncertainty in the coefficients of the underlying model. The first one, known
as Karhunen-Lo\`{e}ve (KL) expansion~\cite{Loeve-1978-PT2}, is suitable for a
finite-dimensional representation of Gaussian random fields, and leads in
conjunction with the stochastic Galerkin finite element discretizations to
block sparse structure of global matrices. Random fields with general
probability distributions are also sometimes represented using the
KL\ expansion, however it might be more suitable to use so-called generalized
polynomial chaos (gPC) expansion~\cite{Xiu-2002-WAP} for their representation.
This in general leads to block dense structure of the global matrices. In
particular, in this paper we will focus on the random coefficients of the
model elliptic PDE with lognormal distribution. Such discretization is done
within the gPC framework, using Hermite polynomials~\cite{Ghanem-1999-NGS}. It
has been shown in~\cite[Theorem 18]{Matthies-2005-GML} that in order to
guarantee a complete representation of the lognormal random field, one has to
use twice the order of polynomial expansion of the coefficient than of the
solution,
which indeed leads to a fully block dense structure of the global stochastic
Galerkin matrix. Clearly, based on the block sparsity pattern of the global
stochastic Galerkin matrix, a certain class of preconditioners might be more
suitable, compared to the others, for a particular problem to be solved.

Probably the most simple, yet quite powerful method is the mean-based
preconditioner proposed by Pellissetti and Ghanem~\cite{Pellissetti-2000-ISS}
and analyzed by Powell and Elman~\cite{Powell-2009-BDP}.
It is a block diagonal preconditioner which uses the information carried by
the mean of the random coefficient and the corresponding stiffness matrix,
i.e., the zeroth order term of the coefficient expansion. Subsequently,
Ullmann has developed a Kronecker product
preconditioner~\cite{Ullmann-2010-KPD}, cf. also~\cite{Ullmann-2012-EIS}, that
makes use of information carried by higher order terms and improves the
convergence quite significantly. The preconditioner has the Kronecker product
structure of the global stochastic Galerkin matrix.
Various iterative methods and preconditioners including multigrid methods,
based on matrix splitting, were compared by Rosseel and Vandewalle
in~\cite{Rosseel-2010-ISS}. Most recently, the authors have developed
in~\cite{Sousedik-2013-HSC} a hierarchical Schur complement preconditioner
that takes advantage of the recursive hierarchy in the structure of the global
stochastic Galerkin matrix.
We note that an interesting approach to solver parallelization has been
proposed by Keese and Matthies~\cite{Keese-2005-HPS}, and a block-triangular
preconditioner for the block sparse case has been proposed recently by Zheng
et al.~\cite{Zheng-2013-BTP}. Finally, we refer to Ernst and
Ullmann~\cite{Ernst-2010-SGM} for a more general study of the stochastic
Galerkin matrices.

In this paper, we propose two novel strategies for constructing
preconditioners for solution of the linear systems with block dense global
matrices to be used with Krylov subspace iterative methods. In particular, we
present a variation on the hierarchical Schur complement preconditioner
developed recently by the authors~\cite{Sousedik-2013-HSC}, and an adaptation
of the symmetric block Gauss-Seidel method. Both preconditioners are built in
a way to account for the recursive hierarchy in the structure of the global
stochastic Galerkin matrix. Because, unlike in the block sparse case, neither
of the submatrices in this hierarchy is block diagonal, we approximate solves
with submatrices by the associated diagonal block solves. This variant yields
versions of the two preconditioners, which will be called as
\emph{approximate}. We note that the approximate versions thus use
\textquotedblleft locally\textquotedblright\ the block diagonal
preconditioning, which allows for a certain level of decoupling. Our numerical
experiments indicate that this approximation might not be necessarily traded
off for slower convergence rates. In the second variation, we take take an
advantage of the decay of the norms of the stiffness matrices obtained from
the polynomial chaos expansion of the coefficients. The decay of the matrices
allows to truncate the matrix-vector multiplications in the action of the
preconditioners, and therefore these versions of the preconditioners will be
called as \emph{truncated}. We also note that one can replace the direct
solvers of the diagonal blocks by iterative ones, possibly with the same
tolerance as for the outer iterations. Doing so, the preconditioners become
variable and one has to make a careful choice of the Krylov iterative method
used for the global (outer) iterations. In general, it is recommended to use
flexible methods such as the flexible conjugate
gradients~\cite{Notay-2000-FCG}, FGMRES~\cite{Saad-1993-FIP}, or
GMRESR~\cite{vanderVorst-1994-GFN}.
Thus neither the global matrix, nor the matrix of the preconditioner need to
be formed explicitly, and we can use the so called MAT-VEC\ operations
introduced in~\cite{Pellissetti-2000-ISS} for all matrix-vector
multiplications. Provided that we have a preconditioner for the diagonal block
solves available, the ingredients of our methods include only the number of
stiffness matrices from the expansion of the random coefficient. Therefore,
the proposed methods can be viewed as minimally intrusive because they can be
built as wrappers around existing solvers for the corresponding deterministic
problem. Finally, we note that choice of a preconditioner for the diagonal
block solves would not change the convergence in terms of outer iterations,
and we will address this topic elsewhere.

The paper is organized as follows. In Section~\ref{sec:model} we introduce the
model problem and its discretization, in Section~\ref{sec:matrices} we discuss
the structure of the stochastic Galerkin matrices, in
Section~\ref{sec:algorithms} we formulate our algorithms, in
Section~\ref{sec:numerical} we present some numerical experiments, and finaly
in Section~\ref{sec:conclusion} we summarize and conclude our work.

\section{Model problem}

\label{sec:model}Let $\left(  \Omega,\mathcal{F},\mu\right)  $ denote the
probability space associated with a physical experiment. We are interested in
the solution of the stochastic linear elliptic boundary value problem, with
stochastic input and deterministic data, given in a bounded domain$~D\subset$
$%
\mathbb{R}
^{d}$, $d=2,3$. The solution is a random function$~u\left(  x,\omega\right)
:\overline{D}\times\Omega\rightarrow%
\mathbb{R}
$ that almost surely (a.s.) satisfies the equation
\begin{align}
-\nabla\cdot\left(  k\left(  x,\omega\right)  \,\nabla u\left(  x,\omega
\right)  \right)   &  =f\left(  x\right)  \qquad\text{in}\ D\times
\Omega,\label{eq:model-1}\\
u\left(  x,\omega\right)   &  =0\qquad\text{on}\ \partial D\times\Omega,
\label{eq:model-2}%
\end{align}
where $f\in L^{2}\left(  D\right)  $, the gradient symbol $\nabla$ denoting
the differentiation with respect to the spatial variables, and%
\[
0<k_{\min}\leq k\left(  x,\omega\right)  \qquad\text{a.s. in }\Omega,\;\forall
x\in\overline{D}.
\]
The function$~k\left(  x,\omega\right)  :$ $\overline{D}\times\Omega
\rightarrow%
\mathbb{R}
$\ is a random scalar field with
a continuous and square-integrable covariance function. We will assume that
the randomness of $k\left(  x,\omega\right)  $ is induced by a set of random
variables $\xi\left(  \omega\right)  =\left\{  \xi\left(  \omega\right)
\right\}  _{i=1}^{N}$\ that are assumed to be independent Gaussian with zero
mean and unit variance. Thus the domain of $k\left(  \cdot,\omega\right)
$\ is restricted to the subset $\left(  \Omega,\mathcal{F}\left(
\mathcal{\xi}\right)  ,\mu\right)  $ of $\left(  \Omega,\mathcal{F}%
,\mu\right)  $ where $\mathcal{F}\left(  \xi\right)  $\ is the $\sigma
-$algebra generated by $\xi$. We will henceforth use $k\left(  x,\xi\right)
$\ instead of \ $k\left(  x,\omega\right)  $, and the solution $u\left(
x,\omega\right)  $ can be also written as $u\left(  x,\xi\right)  $.

In the variational formulation of problem (\ref{eq:model-1})-(\ref{eq:model-2}%
), we consider the solution of the equation
\begin{equation}
u\in U:\ a\left(  u,v\right)  =\left\langle f,v\right\rangle ,\qquad\forall
v\in U. \label{eq:model-variational}%
\end{equation}
where $U$ is a tensor product space defined as
\[
U=H_{0}^{1}\left(  D\right)  \otimes L^{2}\left(  \Omega\right)
,\qquad\left\Vert u\right\Vert _{U}=\sqrt{\mathbb{E}\left[  \int_{D}\left\vert
\nabla u\right\vert ^{2}dx\right]  },
\]
where $\mathbb{E}\left[  \cdot\right]  $ denotes mathematical expectation, and
the bilinear form $a\left(  \cdot,\cdot\right)  $ along with the right-hand
side are defined as
\[
a\left(  u,v\right)  =\mathbb{E}\left[  \int_{D}k\left(  x,\xi\right)
\,\nabla u\cdot\nabla v\,dx\right]  ,\qquad\left\langle f,v\right\rangle
=\mathbb{E}\left[  \int_{D}f\,v\,dx\right]  .
\]
We assume that $k\left(  x,\xi\right)  $\ can be represented, using a
generalized polynomial chaos (gPC) expansion, as%
\begin{equation}
k\left(  x,\xi\right)  =\sum_{i=0}^{M^{\prime}}k_{i}\left(  x\right)  \psi
_{i}\left(  \xi\right)  ,\quad\text{where }k_{i}\left(  x\right)
=\frac{\mathbb{E}\left[  k\left(  x\right)  \psi_{i}\left(  \xi\right)
\right]  }{\mathbb{E}\left[  \psi_{i}^{2}\right]  }, \label{eq:k-expansion}%
\end{equation}
and $\left\{  \psi_{i}\left(  \xi\right)  \right\}  _{i=0}^{M^{\prime}}$ is a
set of $N-$dimensional Hermite polynomials~\cite{Ghanem-1991-SFE}. Similarly,
we will look for an expansion of the solution as
\begin{equation}
u\left(  x,\xi\right)  =\sum_{i=0}^{M}u_{i}\left(  x\right)  \psi_{i}\left(
\xi\right)  , \label{eq:u-expansion}%
\end{equation}
that converges in $L^{2}\left(  \Omega\times D\right)  $ as $M\rightarrow
\infty$. We will, in particular, model $k\left(  x,\xi\right)  $ as a
truncated lognormal process. To this end, let $g\left(  x,\xi\right)
=g_{0}\left(  x\right)  +\sum_{i=1}^{N}\xi_{i}g_{i}\left(  x\right)  $ be a
truncated Karhunen-Lo\`{e}ve expansion of a Gaussian process defined on$~D$
with known mean $g_{0}\left(  x\right)  $ and covariance function$~C_{g}%
\left(  x,y\right)  $, i.e., $g_{i}\left(  x\right)  $\ are weighted
eigenfunctions of
\begin{equation}
\int_{D}C_{g}\left(  x,y\right)  \varphi_{i}\left(  y\right)  \,dy=\lambda
_{i}\varphi_{i}\left(  x\right)  ,\quad\forall x\in D,\qquad\text{and\quad
}g_{i}\left(  x\right)  =\sqrt{\lambda_{i}}\varphi_{i}\left(  x\right)
\text{.} \label{eq:C-integral}%
\end{equation}
Then $k\left(  x,\xi\right)  $ is defined as $k\left(  x,\xi\right)
=\exp\left[  g\left(  x,\xi\right)  \right]  $. The procedure for computing
the coefficients~$k_{i}\left(  x\right)  $ in~(\ref{eq:k-expansion}) has been
derived by Ghanem in~\cite{Ghanem-1999-NGS}. Denoting $\eta_{j}=\xi_{j}-g_{j}%
$, the coefficients are precisely given, cf.~\cite[eq. (33)]{Ghanem-1999-NGS},
as%
\begin{equation}
k_{i}\left(  x\right)  =\frac{\mathbb{E}\left[  \psi_{i}\left(  \eta\right)
\right]  }{\mathbb{E}\left[  \psi_{i}^{2}\right]  }\exp\left[  g_{0}\left(
x\right)  +\frac{1}{2}\sum_{j=1}^{N}\left(  g_{j}\left(  x\right)  \right)
^{2}\right]  .
\end{equation}
According to~\cite{Matthies-2005-GML}, in order for~(\ref{eq:k-expansion}) to
guarantee a complete representation of the lognormal random field, the order
of polynomial expansion of$~k\left(  x,\xi\right)  $ in~(\ref{eq:k-expansion})
should be twice the order of the expansion of the solution$~u\left(
x,\xi\right)  $ in~(\ref{eq:u-expansion}). Denoting the two orders of the
polynomial expansions of $u\left(  x,\xi\right)  $\ and $k\left(
x,\xi\right)  $ as $P$\ and $P^{\prime}$, resp., with $P^{\prime}=2P$, the
total numbers of the gPC\ polynomials are
\begin{equation}
M+1=\frac{\left(  N+P\right)  !}{N!P!},\quad M^{\prime}+1=\frac{\left(
N+P^{\prime}\right)  !}{N!P^{\prime}!}=\frac{\left(  N+2P\right)  !}{N!\left(
2P\right)  !}. \label{eq:M}%
\end{equation}

We will consider approximations to the variational problem
(\ref{eq:model-variational}) given by the finite element discretizations
of$~H_{0}^{1}\left(  D\right)  $. The solution$~u\left(  x,\xi\right)  $ from
(\ref{eq:u-expansion}) will be thus approximated by%
\begin{equation}
u\left(  x,\xi\right)  =\sum_{i=1}^{N_{\text{dof}}}\sum_{j=0}^{M}u_{ij}%
\phi_{i}\left(  x\right)  \psi_{j}\left(  \xi\right)  ,
\label{eq:u-fem-expansion}%
\end{equation}
where $\left\{  \phi_{i}\left(  x\right)  \right\}  _{i=1}^{N_{\text{dof}}}$
is a finite element basis, and $\left\{  \psi_{j}\left(  \xi\right)  \right\}
_{j=0}^{M}$ is the basis of the Hermite polynomials described above.
Substituting the expansions (\ref{eq:k-expansion}) and
(\ref{eq:u-fem-expansion}) into (\ref{eq:model-variational}) yields a
deterministic system of linear equations
\begin{equation}
\sum_{j=0}^{M}\sum_{i=0}^{M^{\prime}}c_{ijk}K_{i}u_{j}=f_{k},\qquad
k=0,\dots,M, \label{eq:global-system}%
\end{equation}
where $(f_{k})_{l}=\mathbb{E}\left[  \int_{D}f\left(  x\right)  \phi
_{l}\left(  x\right)  \psi_{k}\,dx\right]  $, $\left(  K_{i}\right)
_{lm}=\int_{D}k_{i}(x)\phi_{l}(x),\phi_{m}(x)\,dx$, and $c_{ijk}%
=\mathbb{E}\left[  \psi_{i}\psi_{j}\psi_{k}\right]  $. Each one of the
blocks~$K_{i}$ is thus a deterministic\ stiffness matrix given by$~k_{i}%
\left(  x\right)  $ of size $N_{\text{dof}}\times N_{\text{dof}}$, where
$N_{\text{dof}}$ is the number of spatial degrees of freedom. The
system~(\ref{eq:global-system}) is given by a global stochastic Galerkin
matrix of size $\left(  M+1\right)  N_{\text{dof}}\times\left(  M+1\right)
N_{\text{dof}}$, consisting of $N_{\text{dof}}\times N_{\text{dof}}$\ blocks
$K^{\left(  j,k\right)  }$, and it can be written as
\begin{equation}
\left[
\begin{array}
[c]{ccccc}%
K^{\left(  0,0\right)  } & K^{\left(  0,1\right)  } & \cdots &  & K^{\left(
0,M\right)  }\\
& \ddots &  &  & \\
\vdots &  & K^{\left(  k,k\right)  } &  & \vdots\\
&  &  & \ddots & \\
K^{\left(  M,0\right)  } & K^{\left(  M,1\right)  } & \cdots &  & K^{\left(
M,M\right)  }%
\end{array}
\right]  \left[
\begin{array}
[c]{c}%
u_{\left(  0\right)  }\\
\vdots\\
u_{\left(  k\right)  }\\
\vdots\\
u_{\left(  M\right)  }%
\end{array}
\right]  =\left[
\begin{array}
[c]{c}%
f_{\left(  0\right)  }\\
\vdots\\
f_{\left(  k\right)  }\\
\vdots\\
f_{\left(  M\right)  }%
\end{array}
\right]  , \label{eq:global-matrix}%
\end{equation}
where each of the blocks $K^{\left(  j,k\right)  }$ is obtained as
\begin{equation}
K^{\left(  j,k\right)  }=\sum_{i=0}^{M^{\prime}}c_{ijk}K_{i}.
\label{eq:global-matrix-block}%
\end{equation}
We note that the first diagonal block is obtained by the $0-$th order
polynomial chaos expansion and therefore corresponds to the deterministic
problem obtained using the mean value of the coefficient$~k_{0}$. In
particular,
\[
K^{\left(  0,0\right)  }=K_{0}.
\]
In the next section we discuss the structure of the global stochastic Galerkin
matrix from (\ref{eq:global-matrix}) in somewhat more detail.

\section{Hierarchical structure of the matrices}

\label{sec:matrices}
\begin{figure}[ptb]
\centering
\subfigure[{$\ell_t=0$, $$M'=0, nnz = $70$, $\mathrm{n_{MV}}=70$}]
{\includegraphics[width=6.7cm]{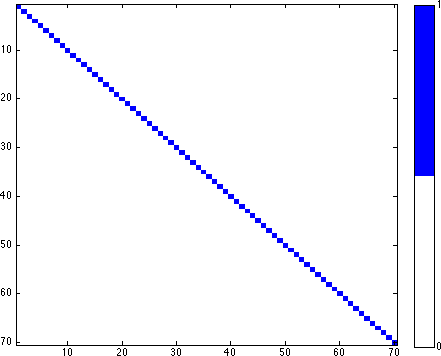} \label{fig:factors-a}}
\subfigure[{$\ell_t=1$, $M'=4$, nnz = $350$, $\mathrm{n_{MV}}=350$}]
{\includegraphics[width=6.7cm]{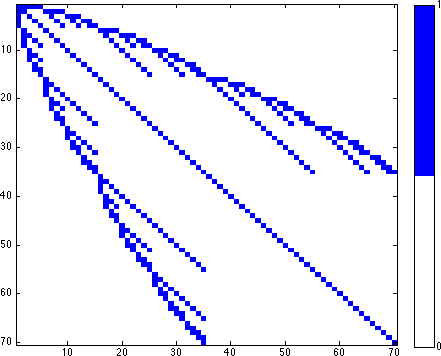} \label{fig:factors-b}}
\newline%
\subfigure[{$\ell_t=2$, $M'=14$, nnz = $1070$, $\mathrm{n_{MV}}=1210$}]
{\includegraphics[width=6.7cm]{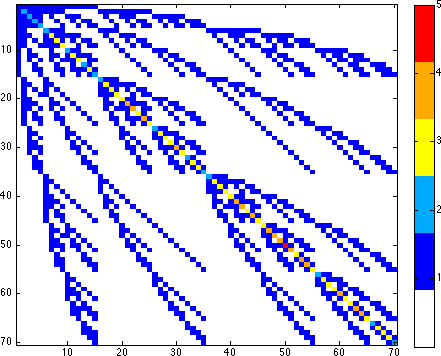} \label{fig:factors-c}}
\subfigure[{$\ell_t=3$, $M'=34$, nnz = $1990$, $\mathrm{n_{MV}}=2610$}]
{\includegraphics[width=6.7cm]{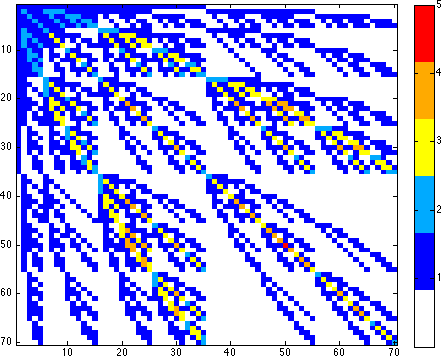} \label{fig:factors-d}}
\newline%
\subfigure[{$\ell_t=4$, $M'=69$, nnz = $3090$, $\mathrm{n_{MV}}=4980$}]
{\includegraphics[width=6.7cm]{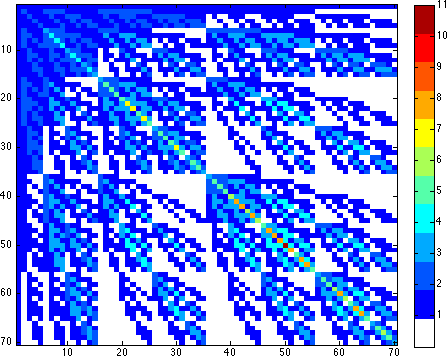} \label{fig:factors-e}}
\subfigure[{$\ell_t=8$, $M'=494$, nnz = $4900$, $\mathrm{n_{MV}}=12,585$}]
{\includegraphics[width=6.7cm]{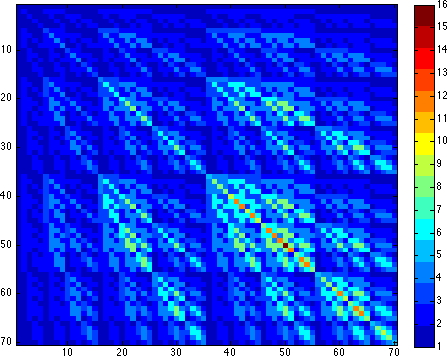} \label{fig:factors-f}}
\caption{Sparsity patterns of the coefficient matrices $c_{P}$ with entries
$c^{(j,k)}=\sum_{i=0}^{M^{\prime}}c_{ijk}$, where $j,k=0,\dots,M$ and
$M+1=\frac{(N+P)!}{N!P!}$. Here $N=P=4$,\ so the size of all the matrices
$c_{P}$ is $M+1=70$, and $M^{\prime}+1=\frac{(N+\ell_{t})!}{N!\ell_{t}!}$,
where $\ell_{t}$ is different in each panel. The colors indicate number of
summations in a position $(j,k)$, $\mathrm{nnz}$ is the number of nonzero
coefficients in the matrix $c_{P}$, and $\mathrm{n_{MV}}$ is the total number
of summations (=MAT-VEC operations).}%
\label{fig:factors}%
\end{figure}

The stochastic Galerkin matrices have recently received
quite a lot of attention, cf.,
e.g.,~\cite{Ernst-2010-SGM,Matthies-2005-GML,Rosseel-2010-ISS}. The key role
in their block structure is played by the constants $c_{ijk}$ and the upper
bound of the summation in equation~(\ref{eq:global-matrix-block}). In general,
there are two types of block sparsity patterns. The first type, typically
regarded as block sparse, is associated with use of only the linear terms
$\xi_{1},\dots,\xi_{N}$ such as appearing in a Karhunen-Lo\`{e}ve expansion of
$k\left(  x,\xi\right)  $, cf. Figure~\ref{fig:factors-b}. The second type,
block dense, is associated with a general (nonlinear in $\xi$'s) form of the
expansion as in~(\ref{eq:k-expansion}). We note that due to our setup with
$P^{\prime}=2P$, we obtain that $M^{\prime}\gg M$ from (\ref{eq:M}), and the
matrices are in general fully dense, cf. Figure~\ref{fig:factors-f}. Recently,
we proposed\ a preconditioner suitable for iterative solution of block
sparse\ stochastic Galerkin systems~\cite{Sousedik-2013-HSC}. Here, we will
focus on a design of preconditioners for more general, block dense, linear systems.

The structure of the global stochastic Galerkin matrix can be understood
through knowledge of the coefficient matrix$~c_{P}$ with entries $c_{P}
\left(  j,k\right)  =\sum_{i=0}^{M^{\prime}}c_{ijk}$ where $j,k=0,\dots,M$,
and the value of$~M$ follows from the fact that we use the $P-$th order
polynomial expansion, cf.~(\ref{eq:M}) \ In general, let us consider
some\ $\ell-$th order polynomial expansion, such that$\ 1\leq\ell\leq P.$\ It
is easy to see that the corresponding coefficient matrix$~c_{\ell}$ will have
a hierarchical structure
\begin{equation}
c_{\ell}=\left[
\begin{array}
[c]{cc}%
c_{\ell-1} & b_{\ell}^{T}\\
b_{\ell} & d_{\ell}%
\end{array}
\right]  ,\qquad\ell=1,\dots,P,
\end{equation}
where $c_{\ell-1}$\ are the first principal submatrices corresponding to the
the $\left(  \ell-1\right)  -$th\ order polynomials expansion. We note that
even though the matrices $c_{\ell}$ are symmetric, the global stochastic
Galerkin matrix in~(\ref{eq:global-matrix}) will be symmetric only if each one
of the stiffness matrices $K_{i}$ is symmetric. Clearly, all matrices $K_{i}$
will have the same sparsity pattern.

\begin{figure}[ptb]
\centering
\subfigure[hierarchical Schur complement preconditioner]
{\includegraphics[width=8cm]{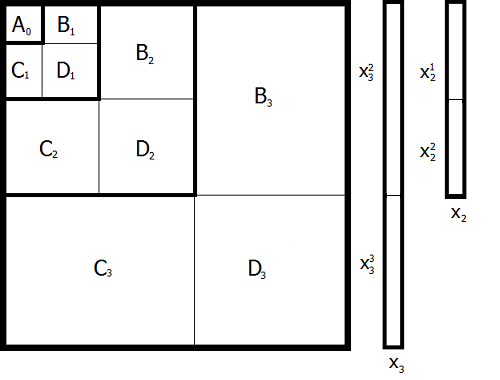} \label{fig:preconditioners-a}}
\subfigure[hierarchical Gauss-Seidel preconditioner]
{\includegraphics[width=8cm]{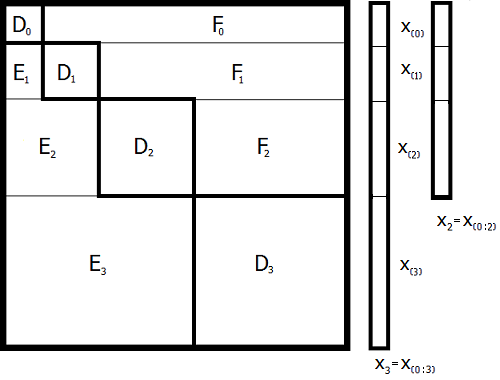} \label{fig:preconditioners-b}}
\caption{Structure of the stochastic Galerkin matrix~(\ref{eq:hierarchy1}) and
of the hierarchical Schur complement preconditioner from
Algorithm~\ref{alg:hS} (left panel), and the structure of the hierarchical
Gauss-Seidel preconditioner from Algorithm~\ref{alg:GS} (right panel). Both
panels also illustrate the structure of vectors introduced in
equation~(\ref{eq:vectors}). }%
\label{fig:preconditioners}%
\end{figure}

In either case, the block sparse or the block dense, the linear system
(\ref{eq:global-matrix}) can be written as
\begin{equation}
A_{P}u_{P}=f_{P}, \label{eq:problemAP}%
\end{equation}
where the global Galerkin matrix $A_{P}$ has the hierarchical structure, cf.
Figure~\ref{fig:preconditioners-a}, given as
\begin{equation}
A_{\ell}=\left[
\begin{array}
[c]{cc}%
A_{\ell-1} & B_{\ell}\\
C_{\ell} & D_{\ell}%
\end{array}
\right]  ,\qquad\ell=P,\dots,1, \label{eq:hierarchy1}%
\end{equation}
and $A_{0}=K_{0}$ is the matrix of the mean.\ Although for the model
problem~(\ref{eq:model-1})-(\ref{eq:model-2}) it holds that $C_{\ell}=B_{\ell
}^{T}$, for$~\ell=1,\dots P$, we will use the general notation
of~(\ref{eq:hierarchy1}) for the sake of generality.

Alternatively, let us consider a hierarchical splitting of the matrix$~A_{P}$,
cf. Figure~\ref{fig:preconditioners-b}, as
\begin{equation}
A_{P}=\left[
\begin{array}
[c]{ccccccc}%
D_{0} &  &  &  & F_{0} &  & \\
& \ddots &  &  &  &  & \\
&  & \ddots &  &  &  & \\
& E_{\ell} &  & D_{\ell} &  & F_{\ell} & \\
&  &  &  & \ddots &  & \\
&  &  &  &  & \ddots & \\
&  & E_{P} &  &  &  & D_{P}%
\end{array}
\right]  . \label{eq:hierarchy2}%
\end{equation}
The blocks $D_{\ell}$, $\ell=1,\dots,P$,\ are the same in (\ref{eq:hierarchy1}%
) and (\ref{eq:hierarchy2}), and we have also, for convenience,\ denoted
$D_{0}=A_{0}$.

The decomposition~(\ref{eq:hierarchy1}) has served as a starting point for the
\emph{hierarchical Schur complement preconditioner}
in~\cite{Sousedik-2013-HSC}, and we will use the decomposition
(\ref{eq:hierarchy2}) to formulate the \emph{hierarchical block symmetric
Gauss-Seidel preconditioner}. Since we are interested here only in the
preconditioners that are \emph{block} and \emph{symmetric} we will drop the
two words for brevity as this shall not cause any confusion. We only note that
in the numerical experiments section we compare a variant of the
\emph{hierarchical} Gauss-Seidel method with the block (non-hierarchical)
Gauss-Seidel method, and the usual (row)\ Gauss-Seidel method is not
considered here at all.

In this paper we make an assumption that it is possible to factorize, e.g. by
the LU-decomposition, the diagonal blocks of the global stochastic Galerkin
matrix or, at least, that we have a preconditioner readily available.


\begin{figure}[ptb]
\centering
\begin{tabular}
[c]{cc}%
\includegraphics[width=7.2cm]{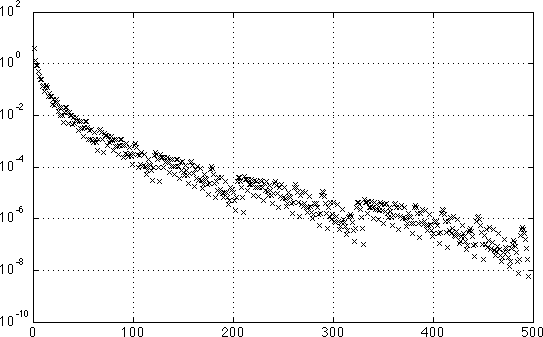} &
\includegraphics[width=7.2cm]{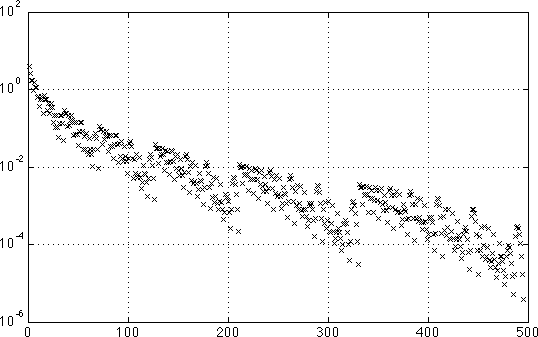}\\
& \\
\includegraphics[width=7.2cm]{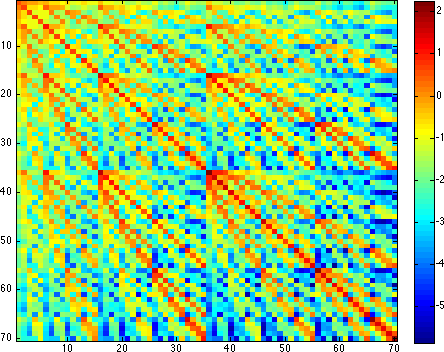} &
\includegraphics[width=7.2cm]{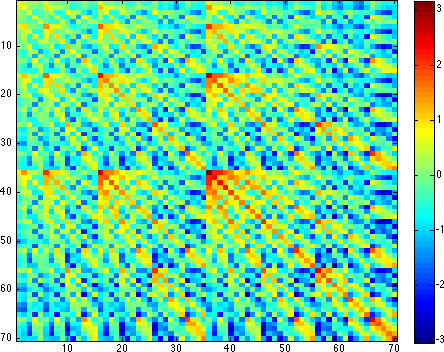}\\
{$CoV=50\%$} & {$CoV=150\%$}%
\end{tabular}
\caption{Norms of the $M^{\prime}+1=495$ stiffness matrices (top panels) in
the representation of the lognormal random field using $N=4$ random variables
and polynomial expansion of order $P^{\prime}=8$, and the decadic logarithm of
the weighted coefficient matrix where in each position $(j,k)$ the coefficient
is obtained as weighted sums $\sum_{i=0}^{494}c_{ijk}\times\mathrm{norm}%
\left(  K_{i}\right)  $ for two coefficients of variation $CoV=50\%$ (left
panel) and $CoV=150\%$ (right panel). The decay of norms of the
matrices~$K_{i}$ illustrated by the top panels motivates the truncation of the
MAT-VEC operations introduced in Algorithm~\ref{alg:tMAT-VEC}, which is used
in the construction of the truncated preconditioners defined in
Algorithms~\ref{alg:hS} and~\ref{alg:GS}. In particular, the construction is
guided by the weighted sums illustrated by the lower panels. }%
\label{fig:norms}%
\end{figure}


\begin{figure}[ptb]
\begin{center}%
\begin{tabular}
[c]{cc}%
\includegraphics[width=6.7cm]{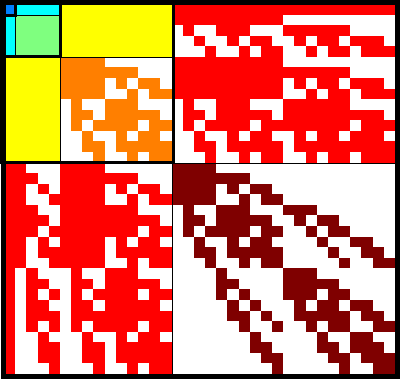} &
\includegraphics[width=6.7cm]{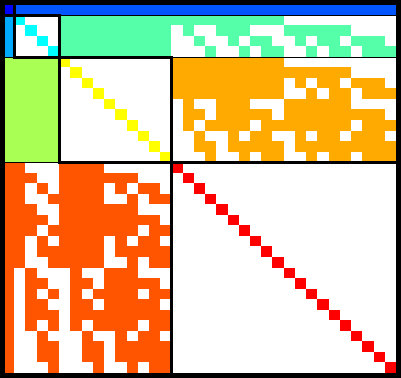}
\end{tabular}
\end{center}
\caption{Block sparsity structure of the hierarchical Schur complement
preconditioner from Algorithm~\ref{alg:hS} (left panel), and the block
sparsity structure of the approximate hierarchical Gauss-Seidel preconditioner
derived from Algorithm~\ref{alg:GS} (right panel). Both preconditioners are
also truncated, setting $\ell_{t}=P$, cf. Figure~\ref{fig:factors-e}. }%
\label{fig:GS}%
\end{figure}


\section{Approximate and truncated hierarchical preconditioners}

\label{sec:algorithms}In the block sparse case, the matrices$~D_{\ell}$ are
block diagonal, as in Figure~\ref{fig:factors-b}. However, considering the
underlying coefficient$~k$ as a random field with lognormal distribution, all
of the matrices$~D_{\ell}$ are block dense, as in Figure~\ref{fig:factors-f}.
The word \emph{approximate} will refer to the fact that we will approximate
the (inverses of the) matrices$~D_{\ell}$ by (inverses of) their diagonal
blocks in the action of our preconditioners. Next, a global matrix-vector
multiplication in an iterative solver can be performed using
formula~(\ref{eq:global-matrix-block}), i.e., one needs to store only the
constants$~c_{ijk}$ and the matrices $K_{i}$ for the so-called called MAT-VEC
operations~\cite{Pellissetti-2000-ISS}. Here, we are interested in
preconditioners that would rely only on the MAT-VEC\ operations as well. The
top panels in Figure~\ref{fig:norms}, show a decaying trend of the norms of
matrices$~K_{i}$ with increasing index$~i$. This indicates, with respect to
the structure of the coefficients $c_{ijk}$, that the lower order terms in the
general $P^{\prime}-$th order expansion might dominate, depending on the
\emph{coefficient of variation} $CoV$, over the higher order terms, cf. lower
panels in Figure~\ref{fig:norms}. Now, the idea is to \emph{truncate}, i.e.,
selectively drop some of the matrices$~K_{i}$ from the MAT-VEC operations in
the action of the preconditioner. There are two possible strategies of this truncation.

First, the \emph{standard truncation }is defined as follows. Let $\ell_{t}$ be
maximal order of the polynomial expansion that we would like to include in the
MAT-VEC\ operations. Then, we define the degree of truncation $M_{t}\leq
M^{\prime}$ as%
\begin{equation}
M_{t}+1=\frac{\left(  N+\ell_{t}\right)  !}{N!\,\ell_{t}!}, \label{eq:M_t}%
\end{equation}
and the matrices $K_{i}$, for all $i=0,\dots,M_{t}$ are included in the MAT-VEC\ operations.

Alternatively, in the \emph{adaptive truncation}\ we select matrices $K_{i}$
that a-priori satisfy certain criterion, e.g., such that their norm is greater
than a given threshold. This is computationaly slightly more demanding because
it requires finding and bookkeeping a subset of matrices$~K_{i}$ where
$i\in\mathcal{M}_{t}\subseteq\left\{  0,\dots,M^{\prime}\right\}  $, however
for higher values of $CoV$, or in the case of a non-monotonous decay of the
matrices $K_{i}$, this might be the preferred truncation strategy.

The algorithm of the truncated block matrix-vector multiplication (tMAT-VEC)
is precisely defined as follows:

\begin{algorithm}
[tMAT-VEC]\label{alg:tMAT-VEC} \ The truncated MAT-VEC\ multiplication
$w=\mathcal{P}_{\ell}v$ is performed as
\begin{equation}
w_{\left(  j\right)  }=\sum_{k=0}^{M}\sum_{i\in\mathcal{M}_{t}}c_{ijk}%
K_{i}v_{\left(  k\right)  },
\end{equation}
where $j,k$ are suitable subsets of $\left\{  0,\dots,M\right\}  $ selecting
blocks of the global Galerkin matrix to be multiplied with, and\ the
truncation is defined by a set $\mathcal{M}_{t}\subseteq\left\{
0,\dots,M^{\prime}\right\}  $ selecting the matrices $K_{i}$ for the MAT-VEC\ operation.
\end{algorithm}

The notation $\mathcal{P}_{\ell}$\ used in Algorithm~\ref{alg:tMAT-VEC} will
correspond in Algorithms~\ref{alg:hS} and~\ref{alg:GS} to either
$\mathcal{B}_{\ell}$, $\mathcal{C}_{\ell}$, $\mathcal{E}_{\ell}$, or
$\mathcal{F}_{\ell}$, and it will denote the corresponding truncated variants
of the block matrix-vector multiplications by the blocks $B_{\ell}$, $C_{\ell
}$, $E_{\ell}$ or $F_{\ell}$ from (\ref{eq:hierarchy1}) and
(\ref{eq:hierarchy2}). We note that the full MAT-VEC\ operation is obtained by
setting $\ell_{t}=P^{\prime}$ in~(\ref{eq:M_t}) so that $M_{t}=M^{\prime}$.


Next, let us introduce some additional notation. Considering
(\ref{eq:hierarchy1}), problem (\ref{eq:problemAP}) can be rewritten as
\begin{equation}
\left[
\begin{array}
[c]{cc}%
A_{P-1} & B_{P}\\
C_{P} & D_{P}%
\end{array}
\right]  \left[
\begin{array}
[c]{c}%
u_{P}^{P-1}\\
u_{P}^{P}%
\end{array}
\right]  =\left[
\begin{array}
[c]{c}%
f_{P}^{P-1}\\
f_{P}^{P}%
\end{array}
\right]  ,
\end{equation}
which motivates the following notation for a vector$~x_{\ell}$, where
$\ell=1,\dots,P$, cf. Figure~\ref{fig:preconditioners}, as%
\begin{equation}
x_{\ell}=\left[
\begin{array}
[c]{c}%
x_{\ell}^{\ell-1}\\
x_{\ell}^{\ell}%
\end{array}
\right]  ,\text{ or }x_{\ell}=\left[
\begin{array}
[c]{c}%
x_{\left(  0\right)  }\\
\vdots\\
x_{\left(  \ell\right)  }%
\end{array}
\right]  \text{.\quad So then }x_{\ell}^{\ell-1}=\left[
\begin{array}
[c]{c}%
x_{\left(  0\right)  }\\
\vdots\\
x_{\left(  \ell-1\right)  }%
\end{array}
\right]  \text{, and }x_{\ell}^{\ell}=x_{\left(  \ell\right)  }\text{.}
\label{eq:vectors}%
\end{equation}
Note that $x_{\left(  0\right)  }=x_{1}^{0}=x_{0}$, and we will also for
brevity denote $x_{\ell}=x_{\left(  0:\ell\right)  }$. We are ready to
formulate now the two preconditioners. First, we recall the algorithm of the
hierarchical Schur complement preconditioner (hS) from~\cite[Algorithm
5]{Sousedik-2013-HSC}. Then, we formulate the algorithm of the hierarchical
Gauss-Seidel preconditioner (hGS).

\begin{algorithm}
[Hierarchical Schur complement preconditioner (hS)]\label{alg:hS} The
preconditioner $M_{hS}:r_{P}\longmapsto v_{P}$\ for system~(\ref{eq:problemAP}%
) is defined as follows:

\noindent\textbf{for} $\ell=P,\ldots1$\textbf{,}

\begin{description}
\item split the residual$~r_{\ell}$, based on the hierarchical structure of
matrices, as
\begin{equation}
r_{\ell}=\left[
\begin{array}
[c]{c}%
r_{\ell}^{\ell-1}\\
r_{\ell}^{\ell}%
\end{array}
\right]  ,
\end{equation}

\item compute the pre-correction as
\begin{equation}
g_{\ell-1}=r_{\ell}^{\ell-1}-\mathcal{B}_{\ell}\mathcal{D}_{\ell}^{-1}r_{\ell
}^{\ell}.
\end{equation}

\item If $\ell>1$, set
\begin{equation}
r_{\ell-1}=g_{\ell-1}.
\end{equation}

\item Else (if $\ell=1)$, solve the system $A_{0}v_{0}=g_{0}$.

\end{description}

\noindent\textbf{end}

\noindent\textbf{for} $\ell=1,\ldots P$\textbf{,}

\begin{description}
\item compute the post-correction, i.e., set $v_{\ell}^{\ell-1}=v_{\ell-1}$,
solve
\begin{equation}
v_{\ell}^{\ell}=\mathcal{D}_{\ell}^{-1}\left(  r_{\ell}^{\ell}-\mathcal{C}%
_{\ell}v_{\ell}^{\ell-1}\right)  ,
\end{equation}

\item and concatenate%
\begin{equation}
v_{\ell}=\left[
\begin{array}
[c]{c}%
v_{\ell}^{\ell-1}\\
v_{\ell}^{\ell}%
\end{array}
\right]  .
\end{equation}

\item If $\ell<P,$ set $v_{\ell+1}^{\ell}=v_{\ell}$.
\end{description}

\noindent\textbf{end}
\end{algorithm}

\begin{algorithm}
[Hierarchical Gauss-Seidel preconditioner (hGS)]\label{alg:GS} The
preconditioner $M_{hGS}:r_{P}\longmapsto v_{P}$\ for
system~(\ref{eq:problemAP}) is defined as follows: \newline set the initial
solution $u_{P}$ and update it in the following steps,
\begin{equation}
v_{\left(  0\right)  }=\mathcal{D}_{0}^{-1}\left(  r_{\left(  0\right)
}-\mathcal{F}_{0}v_{\left(  1:P\right)  }\right)  , \label{eq:algGS1}%
\end{equation}

\noindent\textbf{for} $\ell=1,\ldots P-1$\textbf{,}

\begin{description}
\item
\begin{equation}
v_{\left(  \ell\right)  }=\mathcal{D}_{\ell}^{-1}\left(  r_{\left(
\ell\right)  }-\mathcal{E}_{\ell}v_{\left(  0:\ell-1\right)  }-\mathcal{F}%
_{\ell}v_{\left(  \ell+1:P\right)  }\right)  , \label{eq:algGS2}%
\end{equation}

\end{description}

\noindent\textbf{end}

\noindent\
\begin{equation}
v_{\left(  P\right)  }=\mathcal{D}_{P}^{-1}\left(  r_{\left(  P\right)
}-\mathcal{E}_{P}v_{\left(  0:P-1\right)  }\right)  ,
\end{equation}

\noindent\textbf{for} $\ell=P-1,\ldots1$\textbf{,}

\begin{description}
\item
\begin{equation}
v_{\left(  \ell\right)  }=\mathcal{D}_{\ell}^{-1}\left(  r_{\left(
\ell\right)  }-\mathcal{E}_{\ell}v_{\left(  0:\ell-1\right)  }-\mathcal{F}%
_{\ell}v_{\left(  \ell+1:P\right)  }\right)  ,
\end{equation}

\end{description}

\noindent\textbf{end}

\noindent%
\begin{equation}
v_{\left(  0\right)  }=\mathcal{D}_{0}^{-1}\left(  r_{\left(  0\right)
}-\mathcal{F}_{0}v_{\left(  1:P\right)  }\right)  .
\end{equation}

\end{algorithm}

In our implementation, we initialize $u_{P}=0$ in Algorithm~\ref{alg:GS}. Then
the multiplications of$~\mathcal{F}_{\ell}$, $\ell=0,\dots,P-1$, by this zero
vector will vanish from (\ref{eq:algGS1})-(\ref{eq:algGS2}). This also reduces
the computational cost of the MAT-VEC operations.

We have already mentioned that the \emph{approximate} variants of the two
preconditioners are obtained by replacing the solves with the blocks
$\mathcal{D}_{\ell}$, $\ell=1,\dots,P$, by the corresponding diagonal block
solves, cf. the right panel of Figure~\ref{fig:GS}.

\subsection{Computational cost considerations}

In order to compare the computational cost of the preconditioners, we first
recall two common solution algorithms against which comparisons will be
discussed. First, let us define matrices $G_{\alpha}$, $\alpha=0,\dots
,M^{\prime}$, using the coefficients~$c_{ijk}$ as
\[
G_{\alpha}=c_{\alpha jk},\qquad j=0,\dots,M,\quad k=0,\dots,M.
\]
Let the symbol $\otimes$ denote the Kronecker product. The mean-based
preconditioner~\cite{Pellissetti-2000-ISS,Powell-2009-BDP} is defined as%
\begin{equation}
M_{mb}=\text{diag}\left(  G_{0}\right)  \otimes A_{0}.\label{eq:M_mb}%
\end{equation}
Next, let us define a matrix
\[
G=\sum_{\alpha=0}^{M^{\prime}}\frac{\text{tr}\left(  A_{\alpha}^{T}%
A_{0}\right)  }{\text{tr}\left(  A_{0}^{T}A_{0}\right)  }G_{\alpha}.
\]
The Kronecker product preconditioner~\cite{Ullmann-2010-KPD} is defined as
\begin{equation}
M_{K}=G\otimes A_{0}=\left(  G\otimes I_{N_{\text{dof}}}\right)  \left(
I_{M+1}\otimes A_{0}\right)  =\left(  I_{M+1}\otimes A_{0}^{-1}\right)
\left(  G^{-1}\otimes I_{N_{\text{dof}}}\right)  .\label{eq:M_K}%
\end{equation}
The multiplication of $\left(  I_{M+1}\otimes A_{0}^{-1}\right)  $ by a vector
is essentially the same as an application of the mean-based preconditioner,
and an effficient implementation of $\left(  G^{-1}\otimes I_{\text{ndof}%
}\right)  $\ is described in \cite{MoravitzMartin-2006-SKP}. We can
immediately see from~(\ref{eq:M_mb}) that the mean based-preconditioner
requires $M+1$ solves with the matrix $A_{0}$ of size $\left(  N_{\text{dof}%
}\times N_{\text{dof}}\right)  $, and from~(\ref{eq:M_K}) that the Kronecker
product preconditioner requires in addition $N_{\text{dof}}$ solves with the
matrix$~G$ of size $\left(  \left(  M+1\right)  \times\left(  M+1\right)
\right)  $, which requires $\mathcal{O}\left(  N_{\text{dof}}\times\left(
M+1\right)  ^{2}\right)  $ operations if a Cholesky decomposition of $G$ is
available. Both of our approximate hierarchical preconditioners require
2$\left(  M+1\right)  $ solves with the diagonal blocks of the global Galerkin
matrix, and a certain number of matrix-vector multiplications with stiffness
matrices$~K_{i}$. The number of these multiplications depends on the degree of
truncation of the MAT-VEC\ operations, and it is denoted by $nz\left(
c_{ijk}\right)  $ in Tables~\ref{tab:1} and~\ref{tab:2}. Clearly, the
Kronecker product preconditioner is very efficient if (i)~it is possible to
obtain a Cholesky decomposition of$~G$, and (ii)~the number of spatial degrees
of freedom $N_{\text{dof}}$ is not too large. Moreover, it was assumed
in~\cite{Powell-2009-BDP,Ullmann-2010-KPD} that a solve with $A_{0}$ can be
performed in $\mathcal{O}\left(  N_{\text{dof}}\right)  $\ operations using a
multigrid solver. However, in situations when the multigrid solver is not
suitable or the size of the discretized spatial problem is too large, the
question of selecting block solvers becomes much more delicate. For example,
we might even consider \textquotedblleft inner\textquotedblright\ Krylov
iterations for the solves with the diagonal blocks. The preconditioner would
thus become variable and we would also need to make a careful choice of the
Krylov iterative method used for the global (\textquotedblleft
outer\textquotedblright) iterations. Our initial experiments with such choice
were successful, and the present work is the first step in the development of
such solvers.


\section{Numerical experiments}

\label{sec:numerical}

We have implemented the stochastic Galerkin finite element method for the
model problem~(\ref{eq:model-1})-(\ref{eq:model-2}) on a two-dimensional
domain with dimensions $[0,1]\times\lbrack0,1]$\ uniformly discretized by
$10\times10$\ Lagrangean finite elements. The mean value of the
coefficient$~k$ was set to$~k_{0}=1$, the correlation length $L=0.5$, and we
have used the covariance kernel in~(\ref{eq:C-integral}) as%
\begin{equation}
C_{g}\left(  x,y\right)  =\sigma^{2}\exp\left(  -\frac{\left\Vert
x-y\right\Vert _{1}}{L}\right)  ,
\end{equation}
where $\sigma$ denotes the standard deviation. In the experiments reported in
Tables~\ref{tab:log-N}, \ref{tab:log-P} and~\ref{tab:log-h} we have set
$\sigma=1$, and so the coefficient of variation $CoV=\sigma/k_{0}=100\%$. We
have compared convergence of the flexible conjugate gradients, the mean-based
preconditioner~(mb), the Kronecker product preconditioner (K), the
hierarchical Schur complement preconditioner~(hS), the approximate
hierarchical Schur complement preconditioner~(ahS), symmetric (block)
Gauss-Seidel preconditioner~(GS), and the approximate hierarchical
Gauss-Seidel preconditioner~(ahGS). We note that in all cases we have observed
essentially the same convergence of the standard and flexible conjugate
gradients~\cite{Notay-2000-FCG}.
We have also tested direct and iterative solvers, with the same tolerance as
for the outer iterations, for the inner block solves. Our experiments indicate
that the convergence in terms of outer (global) iteration counts is not
sensitive to the particular choice of inner block solves. The results are
summarized in Tables~\ref{tab:log-N}-\ref{tab:2}. First, we have compared the
convergence using the preconditioners with no truncation, varying either of
the stochastic dimension$~N$, the order of polynomial expansion~$P$, the
coefficient of variation~$CoV$ or the mesh size~$h$, and keeping other
parameters fixed. These results are, respectively, reported in
Tables~\ref{tab:log-N}-\ref{tab:log-h}. Looking at the tables, we see that the
hS preconditioner performs best, followed by the GS preconditioner. On the
other hand, the convergence of the ahS preconditioner quickly deteriorates,
whereas the convergence of the aGS is quite similar to the GS. In fact, it is
quite interesting to note from Table~\ref{tab:log-CoV} that the convergence of
the ahS begins to deteriorate for values of $CoV>25\%$, however the aGS
remains\ comparable to the GS\ preconditioner for values of $CoV$ at
around$\ 100\%$. From Table~\ref{tab:log-h} we see that the dependence on the
mesh size$~h$ is not very significant. Tables~\ref{tab:1} and~\ref{tab:2}
contain results obtained for variable coefficients of variation $CoV=\sigma
_{\log}/\mu_{\log}$ and different truncation strategy of the
MAT-VEC\ operations in the action of the preconditioners, which is guided by
the decay of the norms of the stiffness matrices$~K_{i}$ obtained from the
finite element discretization of the generalized polynomial chaos expansion of
the random coefficient$~k\left(  x,\omega\right)  $, cf.
Figure~\ref{fig:norms}. Table~\ref{tab:1} contains results obtained using the
\emph{standard truncation}. First, we note that setting $\ell_{t}=0$ yields to
use of only the matrix$~K_{0}$ in the action of the ahS, GS, and ahGS
preconditioners. So the resulting preconditioners are in this case block
diagonal and because they are symmetric, their application is the same as
using twice the mean-based preconditioner. On the contrary, the performance of
the hS preconditioner might be slightly better even with $\ell_{t}=0$, in
particular for higher values of $CoV$, but this is because the hS
preconditioner performs solves with the full submatrices$~D_{\ell}$. However,
we can see that this does not correspondingly improve the convergence, and use
of the off-diagonal blocks seems to be warranted. Indeed, the convergence
improves with more included into the MAT-VEC\ operations in the action of all
of the preconditioners. In particular, we see that for values of $CoV$ at
around $25\%$\ it is sufficient to include only the five matrices
corresponding to the linear terms of the coefficient expansion. On the other
hand, even for higher values of $CoV$ there seems to be no reason for the
GS\ and aGS preconditioners\ to use the matrices obtained from the polynomials
of higher order than the expansion of the solution, i.e., the computational
cost of these preconditioners can be reduced more than twice. Finally,
Table~\ref{tab:2} contains results obtained with the \emph{adaptive
truncation}, where the tolerance$~\tau$ has been set such that the matrices
$K_{i}$ for which $\text{max}_{jk}\left(  c_{ijk}\right)  \times
\mathrm{norm}(K_{i})<\tau$ are dropped from the MAT-VEC operations in the
action of the preconditioners. As before, only few matrices need to be used in
order to significantly improve the convergence. Moreover, it appears that
setting the threshold$~\tau$ too low might have a negative impact on the
convergence of the ahS and ahGS\ preconditioners. In particular, rather
surprisingly, even with higher values of $\tau$ the convergence of the
truncated ahGS is comparable to the convergence of GS with no truncation
(setting $\tau=0$). We can in fact see that, e.g., for the values of $CoV$
equal to$\ 100\%$\ and $150\%$, with value $\tau=10$, the convergence of the
ahGS is essentially the same as the convergence of the GS\ preconditioner with
no truncation.

\begin{table}[pth]
\caption{ Convergence of flexible conjugate gradients for the global
matrix~$A$ obtained by the gPC expansion of the lognormal field with
$CoV=\sigma_{\log}/\mu_{\log}=100\%$, and $\mu_{\log}=1$, preconditioned by
the mean-based preconditioner~(mb), Kronecker product preconditioner (K),
hierarchical Schur complement~(hS), approximate hierarchical Schur
complement~(ahS), block Gauss-Seidel~(GS) and approximate hierarchical block
Gauss-Seidel~(ahGS) preconditioners. Polynomial degree is fixed to $P=4$, and
the stochastic dimension $N$ is variable. Here, $ndof$ is the dimension of
$A$, $iter$ is the number of iterations with the relative residual
tolerance~$10^{-8}$, and $\kappa$ is the condition number estimate from the
L\'{a}nczos sequence in conjugate gradients. }%
\label{tab:log-N}%
\centering
\begin{tabular}
[c]{|cc|cc|cc|cc|cc|cc|cc|}\hline
\multicolumn{2}{|c|}{setup} & \multicolumn{2}{|c|}{mb} &
\multicolumn{2}{|c|}{K} & \multicolumn{2}{|c|}{hS} & \multicolumn{2}{|c|}{ahS}
& \multicolumn{2}{|c|}{GS} & \multicolumn{2}{|c|}{ahGS}\\\hline
$N$ & $ndof$ & $it$ & $\kappa$ & $it$ & $\kappa$ & $it$ & $\kappa$ & $it$ &
$\kappa$ & $it$ & $\kappa$ & $it$ & $\kappa$\\\hline
1 & 605 & 48 & 28.76 & 18 & 3.87 & 15 & 3.40 & 15 & 3.40 & 15 & 3.42 & 15 &
3.42\\
2 & 1815 & 61 & 37.16 & 32 & 10.10 & 16 & 3.62 & 27 & 8.06 & 17 & 3.75 & 16 &
3.45\\
3 & 4235 & 62 & 38.07 & 32 & 10.30 & 16 & 3.76 & 31 & 10.77 & 17 & 3.74 & 18 &
4.35\\
4 & 8470 & 66 & 43.65 & 37 & 13.60 & 16 & 4.17 & 38 & 15.28 & 19 & 4.29 & 19 &
4.74\\\hline
\end{tabular}
\end{table}

\begin{table}[pth]
\caption{ Convergence of flexible conjugate gradients for the global
matrix~$A$ preconditioned by the mean-based preconditioner~(mb), Kronecker
product preconditioner (K), hierarchical Schur complement~(hS), approximate
hierarchical Schur complement~(ahS), block Gauss-Seidel~(GS) and approximate
hierarchical block Gauss-Seidel~(ahGS) preconditioners. stochastic dimension
is fixed to $N=4$, and the polynomial degree $P$ is variable. The other
headings are same as in Table~\ref{tab:log-N}. }%
\label{tab:log-P}%
\centering
\begin{tabular}
[c]{|cc|cc|cc|cc|cc|cc|cc|}\hline
\multicolumn{2}{|c|}{setup} & \multicolumn{2}{|c|}{mb} &
\multicolumn{2}{|c|}{K} & \multicolumn{2}{|c|}{hS} & \multicolumn{2}{|c|}{ahS}
& \multicolumn{2}{|c|}{GS} & \multicolumn{2}{|c|}{ahGS}\\\hline
$P$ & $ndof$ & $it$ & $\kappa$ & $it$ & $\kappa$ & $it$ & $\kappa$ & $it$ &
$\kappa$ & $it$ & $\kappa$ & $it$ & $\kappa$\\\hline
1 & 605 & 15 & 3.50 & 14 & 2.44 & 7 & 1.39 & 11 & 1.76 & 8 & 1.39 & 8 & 1.31\\
2 & 18 & 28 & 8.95 & 21 & 4.74 & 10 & 1.93 & 16 & 3.04 & 12 & 1.97 & 11 &
1.76\\
3 & 4235 & 44 & 20.04 & 29 & 8.28 & 13 & 2.80 & 24 & 6.09 & 15 & 2.87 & 14 &
2.58\\
4 & 8470 & 66 & 43.65 & 37 & 13.60 & 16 & 4.17 & 38 & 15.28 & 19 & 4.29 & 19 &
4.74\\\hline
\end{tabular}
\end{table}

\begin{table}[pth]
\caption{ Convergence of flexible conjugate gradients for the global
matrix~$A$ preconditioned by the mean-based preconditioner~(mb), Kronecker
product preconditioner (K), hierarchical Schur complement~(hS), approximate
hierarchical Schur complement~(ahS), block Gauss-Seidel~(GS) and approximate
hierarchical block Gauss-Seidel~(ahGS) preconditioners. Stochastic dimension
and polynomial degree are fixed to $N=P=4$, and the coefficient of variation
$CoV$ is variable $\mu_{\log}=1$. The other headings are same as in
Table~\ref{tab:log-N}. }%
\label{tab:log-CoV}%
\centering
\begin{tabular}
[c]{|c|cc|cc|cc|cc|cc|cc|}\hline
\multicolumn{1}{|c|}{setup} & \multicolumn{2}{|c|}{mb} &
\multicolumn{2}{|c|}{K} & \multicolumn{2}{|c|}{hS} & \multicolumn{2}{|c|}{ahS}
& \multicolumn{2}{|c|}{GS} & \multicolumn{2}{|c|}{ahGS}\\\hline
$CoV\,(\%)$ & $it$ & $\kappa$ & $it$ & $\kappa$ & $it$ & $\kappa$ & $it$ &
$\kappa$ & $it$ & $\kappa$ & $it$ & $\kappa$\\\hline
25 & 16 & 3.24 & 14 & 2.37 & 7 & 1.18 & 8 & 1.25 & 7 & 1.18 & 6 & 1.12\\
50 & 29 & 9.36 & 22 & 4.96 & 10 & 1.78 & 14 & 2.45 & 11 & 1.77 & 10 & 1.62\\
75 & 46 & 22.21 & 30 & 8.73 & 13 & 2.85 & 23 & 6.01 & 15 & 2.82 & 14 & 2.69\\
100 & 66 & 43.65 & 37 & 13.60 & 16 & 4.17 & 38 & 15.28 & 19 & 4.29 & 19 &
4.74\\
125 & 85 & 72.76 & 45 & 19.79 & 19 & 5.54 & 58 & 36.34 & 23 & 5.98 & 26 &
8.21\\
150 & 103 & 107.07 & 52 & 27.02 & 21 & 6.85 & 84 & 77.73 & 26 & 7.75 & 35 &
13.74\\\hline
\end{tabular}
\end{table}

\begin{table}[pth]
\caption{ Convergence of flexible conjugate gradients for the global
matrix~$A$ preconditioned by the mean-based preconditioner~(mb), Kronecker
product preconditioner (K), hierarchical Schur complement~(hS), approximate
hierarchical Schur complement~(ahS), block Gauss-Seidel~(GS) and approximate
hierarchical block Gauss-Seidel~(ahGS) preconditioners. Stochastic dimension
and polynomial degree are fixed to $N=P=4$, the coefficient of variation is
$CoV=100\%$, and the mesh size $h$ is variable. The other headings are same as
in Table~\ref{tab:log-N}. }%
\label{tab:log-h}%
\centering
\begin{tabular}
[c]{|cc|cc|cc|cc|cc|cc|cc|}\hline
\multicolumn{2}{|c|}{setup} & \multicolumn{2}{|c|}{mb} &
\multicolumn{2}{|c|}{K} & \multicolumn{2}{|c|}{hS} & \multicolumn{2}{|c|}{ahS}
& \multicolumn{2}{|c|}{GS} & \multicolumn{2}{|c|}{ahGS}\\\hline
$h$ & $ndof$ & $it$ & $\kappa$ & $it$ & $\kappa$ & $it$ & $\kappa$ & $it$ &
$\kappa$ & $it$ & $\kappa$ & $it$ & $\kappa$\\\hline
$1/5$ & 2520 & 59 & 40.62 & 35 & 14.06 & 15 & 3.84 & 35 & 15.43 & 18 & 3.99 &
19 & 4.91\\
$1/10$ & 8470 & 66 & 43.65 & 37 & 13.60 & 16 & 4.17 & 38 & 15.28 & 19 & 4.29 &
19 & 4.74\\
$1/15$ & 17920 & 68 & 44.42 & 39 & 13.87 & 16 & 4.24 & 39 & 15.81 & 19 &
4.38 & 20 & 4.72\\
$1/20$ & 30870 & 69 & 44.89 & 39 & 14.07 & 17 & 4.25 & 40 & 16.23 & 19 &
4.37 & 20 & 4.78\\
$1/25$ & 47320 & 69 & 44.94 & 40 & 14.13 & 17 & 4.26 & 41 & 16.38 & 20 &
4.40 & 20 & 4.81\\
$1/30$ & 67270 & 71 & 45.11 & 40 & 14.04 & 17 & 4.26 & 41 & 16.38 & 19 &
4.37 & 20 & 4.75\\\hline
\end{tabular}
\end{table}

\begin{table}[ptb]
\caption{Convergence of the flexible conjugate gradients for the global
stochastic Galerkin matrix preconditioned by the mean-based preconditioner
(mb), Kronecker product preconditioner (K), \ hierarchical Schur complement
(hS), approximate hierarchical Schur complement (ahS), block Gauss-Seidel (GS)
and approximate hierarchical block Gauss-Seidel (ahGS) preconditioners with a
variable degree of truncation of the MAT-VEC operations in the action of the
preconditioner. Here $\ell_{t}$ is the maximum polynomial order of the
coefficient expansion used in the construction of the preconditioner, so that
$M_{t}+1$ is the degree of truncation of the MAT-VEC operations i.e., the
number of retained matrices, $nz(c_{ijk})$ is the number of nonzeros in the
truncated tensor~$c_{ijk}$, $iter$ is the number of iterations with the
relative residual tolerance~$10^{-8}$, and $\kappa$ is the condition number
estimate from the L\'{a}nczos sequence in flexible conjugate gradients.}%
\label{tab:1}
\begin{center}%
\begin{tabular}
[c]{|ccc|cc|cc|cc|cc|}\hline
\multicolumn{3}{|c|}{setup} & \multicolumn{2}{|c|}{hS} &
\multicolumn{2}{|c|}{ahS} & \multicolumn{2}{|c|}{GS} &
\multicolumn{2}{|c|}{ahGS}\\\hline
$\ell_{t}$ & $M_{t}+1$ & $nz(c_{ijk})$ & $it$ & $\kappa$ & $it$ & $\kappa$ &
$it$ & $\kappa$ & $it$ & $\kappa$\\\hline
\multicolumn{11}{|c|}{$CoV=25\%$ \quad( mb: \, $it=16$ \quad$\kappa=3.24$,
\,\, K: \, $it=14$ \quad$\kappa=2.37$ )}\\\hline
0 & 1 & 70 & 16 & 3.20 & 16 & 3.19 & 16 & 3.19 & 16 & 3.19\\
1 & 5 & 350 & 8 & 1.27 & 8 & 1.33 & 7 & 1.23 & 7 & 1.23\\
2 & 15 & 1210 & 7 & 1.21 & 8 & 1.25 & 7 & 1.20 & 6 & 1.17\\
3 & 35 & 2610 & 7 & 1.18 & 8 & 1.25 & 7 & 1.17 & 6 & 1.11\\
4 & 70 & 4980 & 7 & 1.18 & 8 & 1.25 & 7 & 1.18 & 6 & 1.12\\
8 & 495 & 12585 & 7 & 1.18 & 8 & 1.25 & 7 & 1.18 & 6 & 1.12\\\hline
\multicolumn{11}{|c|}{$CoV=50\%$ \quad( mb: \, $it=29$ \quad$\kappa=9.36$,
\,\, K: $it=22$ \quad$\kappa=4.96$)}\\\hline
0 & 1 & 70 & 28 & 8.90 & 28 & 8.99 & 28 & 8.99 & 28 & 8.99\\
1 & 5 & 350 & 14 & 2.54 & 15 & 2.83 & 14 & 2.54 & 14 & 2.54\\
2 & 15 & 1210 & 13 & 2.47 & 14 & 2.74 & 12 & 2.15 & 11 & 2.12\\
3 & 35 & 2610 & 10 & 1.79 & 14 & 2.44 & 10 & 1.68 & 10 & 1.69\\
4 & 70 & 4980 & 10 & 1.79 & 13 & 2.39 & 11 & 1.78 & 10 & 1.68\\
8 & 495 & 12585 & 10 & 1.78 & 14 & 2.45 & 11 & 1.77 & 10 & 1.62\\\hline
\multicolumn{11}{|c|}{$CoV=100\%$ \quad( mb: \, $it=66$ \quad$\kappa=43.65$,
\,\, K: $it=37$ \quad$\kappa=13.60$)}\\\hline
0 & 1 & 70 & 53 & 32.88 & 58 & 39.94 & 58 & 39.94 & 58 & 39.94\\
1 & 5 & 350 & 32 & 11.83 & 35 & 14.29 & 33 & 12.88 & 33 & 12.88\\
2 & 15 & 1210 & 32 & 12.31 & 35 & 15.13 & 28 & 9.35 & 29 & 10.05\\
3 & 35 & 2610 & 20 & 5.03 & 36 & 14.32 & 19 & 4.64 & 22 & 6.23\\
4 & 70 & 4980 & 20 & 5.24 & 32 & 11.92 & 19 & 4.55 & 20 & 5.33\\
8 & 495 & 12585 & 16 & 4.17 & 38 & 15.28 & 19 & 4.29 & 19 & 4.74\\\hline
\multicolumn{11}{|c|}{$CoV=150\%$ \quad( mb: \, $it=103$ \quad$\kappa
=107.067$, \,\, K: $it=52$ \quad$\kappa=27.0203$)}\\\hline
0 & 1 & 70 & 71 & 61.44 & 89 & 90.18 & 89 & 90.18 & 89 & 90.18\\
1 & 5 & 350 & 51 & 29.92 & 59 & 39.66 & 57 & 36.85 & 57 & 36.85\\
2 & 15 & 1210 & 51 & 30.18 & 60 & 42.06 & 46 & 24.26 & 50 & 27.53\\
3 & 35 & 2610 & 31 & 11.06 & 71 & 55.98 & 30 & 10.17 & 40 & 19.53\\
4 & 70 & 4980 & 32 & 12.05 & 58 & 38.08 & 28 & 9.42 & 34 & 13.86\\
8 & 495 & 12585 & 21 & 6.85 & 84 & 77.73 & 26 & 7.75 & 35 & 13.74\\\hline
\end{tabular}
\end{center}
\end{table}


\begin{table}[ptb]
\caption{Convergence of the flexible conjugate gradients for the global
stochastic Galerkin matrix preconditioned by the mean-based preconditioner
(mb), Kronecker product preconditioner (K), hierarchical Schur complement
(hS), approximate hierarchical Schur complement (ahS), block Gauss-Seidel (GS)
and approximate hierarchical block Gauss-Seidel (ahGS) preconditioners with a
variable degree of truncation of the MAT-VEC operations in the action of the
preconditioner. Here $\tau$ is the tolerance such that the matrices $K_{i}$
for which $\text{max}_{jk}\left(  c_{ijk}\right)  \times\mathrm{norm}%
(K_{i})<\tau$ are dropped from the MAT-VEC operations in the action of the
preconditioner and $N_{\text{adapt}}$ is the number of retained matrices,
$nz(c_{ijk})$ is the number of nonzeros in the truncated tensor $c_{ijk}$,
$iter$ is the number of iterations with the relative residual
tolerance~$10^{-8}$, and $\kappa$~is the condition number estimate from the
L\'{a}nczos sequence in flexible conjugate gradients.}%
\label{tab:2}
\begin{center}%
\begin{tabular}
[c]{|ccc|cc|cc|cc|cc|}\hline
\multicolumn{3}{|c|}{setup} & \multicolumn{2}{|c|}{hS} &
\multicolumn{2}{|c|}{ahS} & \multicolumn{2}{|c|}{GS} &
\multicolumn{2}{|c|}{ahGS}\\\hline
$\tau$ & $N_{\text{adapt}}$ & $nz(c_{ijk})$ & $it$ & $\kappa$ & $it$ &
$\kappa$ & $it$ & $\kappa$ & $it$ & $\kappa$\\\hline
\multicolumn{11}{|c|}{$CoV = 25\%$ \qquad mb: \, $it=16$ \, $\kappa=3.24$
\quad K: \, $it=14$ \, $\kappa=2.37$}\\\hline
$10$ & 5 & 345 & 10 & 1.63 & 10 & 1.58 & 10 & 1.57 & 10 & 1.57\\
$1$ & 13 & 877 & 7 & 1.20 & 8 & 1.24 & 7 & 1.18 & 6 & 1.12\\
$0.1$ & 32 & 2057 & 7 & 1.18 & 8 & 1.24 & 7 & 1.17 & 6 & 1.12\\
$0$ & 495 & 12585 & 7 & 1.18 & 8 & 1.25 & 7 & 1.18 & 6 & 1.12\\\hline
\multicolumn{11}{|c|}{$CoV = 50\%$ \qquad mb: \, $it=29$ \quad$\kappa=9.36$
\quad K: \, $it=22$ \, $\kappa=4.96$}\\\hline
$10$ & 11 & 677 & 12 & 2.18 & 13 & 2.24 & 11 & 1.83 & 11 & 1.83\\
$1$ & 32 & 1958 & 11 & 1.91 & 13 & 2.25 & 11 & 1.77 & 10 & 1.62\\
$0.1$ & 86 & 4765 & 10 & 1.78 & 14 & 2.41 & 11 & 1.77 & 10 & 1.62\\
$0$ & 495 & 12585 & 10 & 1.78 & 14 & 2.45 & 11 & 1.77 & 10 & 1.62\\\hline
\multicolumn{11}{|c|}{$CoV = 100\%$ \qquad mb: \, $it=66$ \quad$\kappa=43.65$
\quad K: \, $it=37$ \, $\kappa=13.60$}\\\hline
$100$ & 2 & 92 & 54 & 34.95 & 63 & 47.94 & 65 & 51.28 & 65 & 51.28\\
$10$ & 25 & 1241 & 27 & 9.12 & 25 & 7.14 & 23 & 6.98 & 19 & 4.56\\
$1$ & 97 & 4731 & 18 & 4.48 & 33 & 12.29 & 19 & 4.55 & 18 & 4.10\\
$0.1$ & 219 & 8202 & 17 & 4.18 & 37 & 14.90 & 19 & 4.31 & 19 & 4.66\\
$0$ & 495 & 12585 & 16 & 4.17 & 38 & 15.28 & 19 & 4.29 & 19 & 4.74\\\hline
\multicolumn{11}{|c|}{$CoV = 150\%$ \qquad mb: \, $it=103$ \, $\kappa=107.07$
\quad K: \, $it=52$ \, $\kappa=27.02$}\\\hline
$100$ & 10 & 336 & 66 & 50.22 & 68 & 48.84 & 70 & 50.68 & 70 & 50.68\\
$10$ & 55 & 2450 & 30 & 10.89 & 44 & 20.38 & 29 & 9.49 & 25 & 6.95\\
$1$ & 171 & 6338 & 23 & 7.00 & 73 & 59.57 & 27 & 7.98 & 32 & 11.48\\
$0.1$ & 313 & 9714 & 22 & 6.86 & 83 & 76.38 & 26 & 7.74 & 35 & 13.54\\
$0.01$ & 436 & 11741 & 21 & 6.84 & 84 & 77.69 & 26 & 7.75 & 35 & 13.72\\
$0$ & 495 & 12585 & 21 & 6.85 & 84 & 77.73 & 26 & 7.75 & 35 & 13.74\\\hline
\end{tabular}
\end{center}
\end{table}

\section{Conclusion}

\label{sec:conclusion}We have proposed two novel strategies for constructing
preconditioners for the iterative solution of the systems of linear algebraic
obtained from the stochastic Galerkin finite element discretizations. Our main
focus was on a class of problems with coefficients characterized by arbitrary
distributions that generally yield\ fully block dense structure of global
stochastic Galerkin matrices. The preconditioners take an advantage of the
hierarchical structure of these matrices. We have, in particular, examined
variants of the hierarchical Schur complement preconditioner proposed recently
by the authors~\cite{Sousedik-2013-HSC}, and variants of the hierarchical
block symmetric Gauss-Seidel preconditioner.

The first variant, called an \emph{approximate}, replaces solves with
submatrices by the associated diagonal block solves. This variant thus
combines global Gauss-Seidel method with \textquotedblleft
local\textquotedblright\ block-diagonal preconditioner, which allows
decoupling of blocks, and introduces a possibility for a parallelisation in an
implementation. Numerical experiments with our model problem indicate that
whereas the performance of the approximate version of the hierarchical Schur
complement preconditioner deteriorates with the increasing values of the
coefficient of variation, the convergence of the approximate hierarchical
Gauss-Seidel preconditioner is quite comparable to the (non-hierarchical)
block Gauss-Seidel preconditioner for the values of the coefficient of
variation equal up to $100\%.$

The second variant, called \emph{truncated}, is based on a truncation of the
sequence of block matrix-vector multiplications, called MAT-VEC\ operations,
used in the action of the preconditioners. The truncation can be performed
using either a standard or adaptive strategy, based on the monotonicity in the
decay of the stiffness matrices, and further alleviates the computational cost
of the preconditioners. Our numerical experiments indicate that truncation, in
particular with adaptive strategy, might not be necessarily traded off for
slower convergence rates.

We have therefore proposed two strategies that are combined for optimal
performance. One strategy introduces decoupled diagonal block solves, and the
other reduces the overall computational cost associated with the action of a
preconditioner. Considering that a multiplication of a vector by the global
stochastic Galerkin matrix in each iteration of a Krylov subspace method is,
of course, performed by the full MAT-VEC\ multiplication with no truncation,
in particular the truncated version of the approximate hierarchical
Gauss-Seidel preconditioner preconditioners thus offers an appealing
combination of good convergence rates and a reasonable computational cost.

\acknowledgements Support from DOE/ASCR is gratefully acknowledged. B.
Soused\'{\i}k has been also supported in part by the Grant Agency of the Czech
Republic \mbox{GA \v{C}R 106/08/0403}.

\bibliographystyle{IJ4UQ_Bibliography_Style}
\bibliography{truncated}

\end{document}